\documentclass[letterpaper, 10 pt, conference]{ieeeconf}  

\IEEEoverridecommandlockouts                              
\usepackage{amsmath} 
\usepackage{amssymb}  
\usepackage{amsfonts}
\usepackage{graphicx}

\title{\LARGE \bf
Convex Optimization In Identification Of\\ 
Stable Non-Linear State Space Models
}

\IEEEoverridecommandlockouts
\author{
Mark M. Tobenkin$^1$ \ Ian R. Manchester$^1$ \ Jennifer Wang$^2$ \ Alexandre Megretski$^3$ \  Russ Tedrake$^1$\\  \ \\
1: Computer Science and Artificial Intelligence Laboratory, Massachusetts Inst. of Tech.\\ 
2:  Department of Brain and Cognitive Sciences, Massachusetts Inst. of Tech.\\
3:  Laboratory for Information and Decision Systems, Massachusetts Inst. of Tech.\\
\{mmt, irm, jenwang, ameg, russt\}@mit.edu
\thanks{This was supported by National Science Foundation Grant No. 0835947.}
}

\newtheorem{thm}{Theorem}

\newtheorem{lem}{Lemma}
\newcounter{remark}
\newcounter{example}

\newenvironment{pf}{\smallbreak\noindent{\it Proof. }}{\hfill$\Box$\smallbreak}
\newenvironment{pf*}[1]{\smallbreak\noindent{\it #1}}{\hfill$\Box$\smallbreak}

\makeatother
\ifx\Box\undefined
  \makeatletter\input{oldlfont.sty} \makeatother
\fi

\newcommand{\rf}[1]{(\ref{#1})}
\newcommand{\cl}[1]{{\cal #1}}

\def\Box{\hbox{\hskip 1pt \vrule width 8pt height 6pt depth 1.5pt
  \hskip 1pt}}

\newcommand{\e}{\epsilon}

\renewcommand{\d}{\delta}
\newcommand{\D}{\Delta}

\renewcommand{\t}{\theta}

\newcommand{\s}{\sigma}

\newcommand{ \R}{\mathbb R}

\newcommand{\EQ}[2]{\begin{equation}\label{#1}#2\end{equation}}

\DeclareRobustCommand{\greektext}{%
  \fontencoding{LGR}\selectfont\def\encodingdefault{LGR}}
\DeclareRobustCommand{\textgreek}[1]{\leavevmode{\greektext #1}}
\DeclareFontEncoding{LGR}{}{}

\newcommand{\RR}{\mathbb{R}}

\usepackage{verbatim}

\begin{document}

\maketitle
\thispagestyle{empty}
\pagestyle{empty}

\begin{abstract}
A new framework for nonlinear system identification is presented in
terms of optimal fitting of stable nonlinear state space equations to
input/output/state data, with a performance objective defined as a
measure of robustness of the simulation error with respect to equation
errors.  Basic definitions and analytical results are presented. The
utility of the method is illustrated on a simple simulation example as
well as experimental recordings from a live neuron.
\end{abstract}

\section{Introduction}
Converting numerical data, originating from either 
physical measurements or computer simulations, to compact
mathematical models is a common challenge in engineering.
The case of {\sl static system identification}, where 
models $y=h(u)$ defined by ``simple" functions
$h(\cdot)$
are fitted to data records of $u$ and $y$,
is a major topic of research in {\sl statistics} and
{\sl machine learning}. This paper is focused on
a subset of
{\sl dynamic system identification} tasks, where 
state space models of the form
\begin{eqnarray}
x(t+1)&=&f(x(t),u(t)),\label{mi1}\\
y(t)&=&g(x(t),u(t)),\label{mi2}
\end{eqnarray}
where $f$ and $g$ are ``simple" functions,
are extracted from data records of $\tilde x, \tilde y,$ and $\tilde u$. We will also consider
continuous time models of the form
\begin{eqnarray}
\dot x(t)&=&f(x(t),u(t)),\label{mi3}
\end{eqnarray}
with \eqref{mi2}. If the data come from simulations of a complex
model rather than experiments then this task is referred to as {\em model reduction}.

There are two common and straightforward approaches to this problem:
\begin{enumerate}
\item {\sl Equation-error minimization:} (see, e.g., \cite{LjungBook,
    Sjoberg95}) A model is sought which minimizes a
  cost function of the following form
\[
\sum_t|\tilde x(t+1)-f(\tilde x(t),\tilde u(t))|^2,
\]
or similar, over the unknown parameters of $f(\cdot)$. A similar optimization can be
set up for $g(\cdot)$. This is typically very cheap computationally:
 if $f(\cdot)$ and $g(\cdot)$ are linear in the unknown parameters then the problem
 reduces to basic
 least squares. However, if there is no incremental stability
 requirement then small equation errors $|\tilde x(t+1)-f(\tilde x(t),\tilde u(t))|$ do not imply small simulation
 errors over extended time intervals. For large scale and nonlinear problems it
 is not unusual to find unstable models by this method, particularly
 if the true system is not in the model class being searched over.
\item {\sl Simulation-error minimization:} (see, e.g., \cite{Farina10}) One sets up a nonlinear
  programming problem to find
\[
\min_\eta\sum_t|\tilde y(t)-y_\eta(t)|^2
\]
where $y_\eta$ is the output of the simulation of the model system
with a particular set of parameters
$\eta$ which define $f(\cdot)$ and $g(\cdot)$. If successful, this can give a more robust fit than
  equation error. However, even if the system equations $f(\cdot), g(\cdot)$ are linear in
  the unknown parameters, the relationship between the unknown
  parameters $\eta$ and the long-term simulation $y_\eta(t)$ will be highly
  nonlinear and the optimization nonconvex. For systems with a large
  number of parameters, this can make global optimization of simulation error
  very difficult unless good
  initial parameter guesses are known, which is seldom the case when
  considering black box model structures.
\end{enumerate}
The method proposed in this paper can be considered a middle-ground
between these two extremes: we formulate a convex optimization
problem to minimize an {\em upper bound} on the true simulation error while
guaranteeing the stability and well-posedness of the identified model. Furthermore, we show
that that in some simple cases the upper bound is tight.



While ensuring stability complicates identification
of both linear and nonlinear models, it is most challenging
in the nonlinear case. Some recently proposed methods include LMI
conditions for linear systems  \cite{Lacy03}, convex relaxations for
linear \cite{Meg} and Wiener-Hammerstein systems \cite{Sou08}, as well
as passivity-like
conditions for linear \cite{KMD} and nonlinear models \cite{Meg08}.

We do not in general guarantee finding statistically optimal or
unbiased estimates. However, for nonlinear or high-order linear systems stability
of the model and
reduction of the long-term error dynamics are often major problems; these have
been our primary targets. 


\section{Problem Setup}
\subsection{State Space Models}
We examine discrete time (DT) state space models
of the form
\begin{eqnarray}
e(x(t+1))&=&f(x(t),u(t)),\label{mi5}\\
y(t)&=&g(x(t),u(t)),\label{mi6}
\end{eqnarray}
where $e:\ \RR^n\mapsto\RR^n$, 
$f:\ \RR^n\times\RR^m\mapsto\RR^n$,
and $g:\ \RR^n\times\RR^m\mapsto\RR^k$
are continuously differentiable functions
such that the equation $e(z)=w$ has a unique
solution $z\in\RR^n$ for every $w\in\RR^n$.

\subsection{Stability}
We consider the DT model \rf{mi5},\rf{mi6} {\sl stable}
if the difference $\{y_1(t)-y_0(t)\}_{t=1}^\infty$ is
square summable for every two solutions 
$(u,x,y)=(u_1,x_1,y_1)$
and $(u,x,y)=(u_0,x_0,y_0)$ of \rf{mi5},\rf{mi6} with
the same input $u_1=u_0=u$. 
This definition can be qualified as that of {\sl 
global incremental output $\ell_2$-stability}.

\subsection{Data}
In applications, we expect to have
input/state/output information available in the form of 
{\sl sampled data}
$\tilde {\cl Z}=\{\tilde z(t_i)\}_{i=1}^N$, where $\tilde z(t_i)=(\tilde v(t_i),\tilde x(t_i),\tilde u(t_i),\tilde y(t_i))$. Here
$\tilde x,\tilde u,\tilde y$ represent approximated
samples of state, input, and output respectively.
Section \ref{sec:state} will discuss approaches for approximating the
state of the system from input-output data.

For the purpose of theoretical
analysis, we will assume that the
input/state/output information is available in the form of 
{\sl signal data}
where $\tilde v$, $\tilde x$,
$\tilde u$, and $\tilde y$ are {\sl signals}, i.e.
functions
\EQ{mi44}{
\tilde x,\tilde v:\ \cl T\mapsto\RR^n,\ 
\tilde u:\ \cl T\mapsto\RR^m,\ 
\tilde y:\ \cl T\mapsto\RR^k
}
such that 
\EQ{mi45}{\cl T=\{1,\dots,N\},\ \ \ 
\tilde v(t-1)=\tilde x(t)\ \forall\ t\in\{2,\dots,N\}.
}

 \subsection{Simulation Error}
 Given DT {\sl signal data}
 $\tilde{\cl Z}$
 and functions $e,f,g$, the {\sl simulation error}
 associated with a model matching \rf{mi5},\rf{mi6} is defined
 as 
 \EQ{mi21}{ \bar{\cl E}=\sum_{t=1}^N|\tilde y(t)- y(t)|^2,}
 where $y(t)$ is defined by \rf{mi5},\rf{mi6} with
 $u(t)\equiv\tilde u(t)$ and
 $x(1)=\tilde x(1)$.

 \subsection{Linearized Simulation Error}
 \label{sec:linearizeddt}
A simplified version of the simulation error measure
$\bar{\cl E}$ is the {\sl linearized simulation error}
$\bar{\cl E}^0$ defined in the following way. 
Consider a ``perturbed'' version of the system equations
\begin{eqnarray}
e(x_\t(t+1))&=&f(x_\t(t),u(t))-f_0(t),\label{mi31}\\
y_\t(t)&=&g(x_\t(t),u(t))-g_0(t).\label{mi32}
\end{eqnarray}

Here $\t \in [0,1]$, $f_0(t)=(1-\t)\e_x(\tilde z(t))$,
and $g_0(t)=(1-\t)\e_y(\tilde z(t))$, where $\epsilon_x$ and $\epsilon_y$ are
the {\sl equation errors} are defined by:
\begin{eqnarray}
\e_x(\tilde z)&=&f(\tilde x,\tilde u)-e(\tilde v),\ \ \ \ \label{mi8}\\ 
\e_y(\tilde z)&=&g(\tilde x,\tilde u)-\tilde y.\label{mi10}
\end{eqnarray}
We examine the solution $( x_\t, y_\t)$ of
\rf{mi31},\rf{mi32} with $x_\t(1) = \tilde x(1)$, $u(t) \equiv \tilde
u(t)$.

By construction, $y_\t=y$ for $\t=1$, 
and $y_\t=\tilde y$ for $\t=0$.
We define
\EQ{mi34}{ \bar{\cl E}^0(\tilde z(t))=\lim_{\t\to0}\frac1{\t^2}
\sum_{t=1}^N|\tilde y(t)- y_\t(t)|^2}
to quantify local sensitivity of model
equations with respect to equation errors.

Using standard linearization analysis, it is easy to
produce alternative expressions for $\bar{\cl E}^0$:
\EQ{mi36}{ \bar{\cl E}^0=
\sum_{t=1}^N|G(\tilde x(t),\tilde u(t))\tilde\D(t)+
\e_y(\tilde z(t))|^2,}
where $\tilde\D(\cdot)$ is defined by
\EQ{mi37}{ E(\tilde x(t+1))\tilde\D(t+1)=
F(\tilde x(t),\tilde u(t))\tilde\D(t)+
\e_x(\tilde z(t)),}
with initial condition $\tilde\D(1)=0$, and
$E = E(x),F=F(x,u)$ and $G=G(x,u)$ defined to be the 
Jacobians (with respect to $x$)
of $e, f$ and $g$ respectively.

\subsection{Optimization Setup}
Within the framework of this paper,
we consider efficient global minimization 
of the simulation error
$\bar{\cl E}$ 
(over all
model functions $e,f,g$, defining a stable system)
as an ultimate (if perhaps unattainable) goal.
We proceed by defining upper bounds
for $\bar{\cl E}$ and $\bar{\cl E}^0$ 
which can be minimized efficiently
by means of convex optimization 
(semidefinite programming). We will also prove 
some theoretical
statements certifying quality of these upper bounds.


\section{Robust Identification Error}
The dependence of the simulation error 
($\bar{\cl E}$ or $\bar{\cl E}^0$) on the coefficients of
system equations \rf{mi5},\rf{mi6} is 
complicated enough to make it a challenging object for
efficient global minimization, especially
under the stability constraint. The objective of this
section is to introduce several versions of
{\sl robust identification error} (RIE) -
a sample-wise measure
of simulation error, motivated by the idea of
using storage functions
and dissipation inequalities to generate
useful upper bounds of $\bar{\cl E}$ and $\bar{\cl E}^0$.


\subsection{Global RIE}
The global RIE measure for a DT model \rf{mi5},\rf{mi6} 
is a function of the coefficients of \rf{mi5},\rf{mi6},
a single data sample 
 \EQ{mi41}{
 \tilde z=(\tilde v,\tilde x,\tilde u,\tilde y)\in\RR^n\times\RR^n\times\RR^m\times\RR^k,
 }
 and an auxiliary parameter $Q=Q'>0$, 
 a positive definite symmetric
 $n$-by-$n$ matrix (for convenience, 
 we only indicate the dependence on $\tilde z$ and $Q$):
\EQ{mi13}{
\cl E_Q(\tilde z)=
\sup_{\D}\left\{|f(\tilde x+\D,\tilde u)-e(\tilde v)|^2_Q-
|\d_e|^2_Q+|\d_y|^2\right\}.
}
where $|a|^2_Q$ is a shortcut for $a'Qa$, and
\EQ{mi42}{ \d_y=g(\tilde x+\D,\tilde u)-\tilde y,\ \ \d_e=e(\tilde x+\D)-e(\tilde x).}
 The following statement explains the
 utility of the RIE measure in generating upper bounds
 of simulation error.

\begin{thm}\label{thm:riegdt}
{\sl The inequality
\EQ{mi49}{  \bar{\cl E}\le\sum_{t=1}^N
\cl E_Q(\tilde z(t)),
}
holds
for every $Q=Q'>0$ and signal data \rf{mi44},\rf{mi45}.
}
\end{thm}
\begin{pf}
By the definition of $\cl E_Q(\tilde z(t))$ we have
  \begin{flalign}
    |f(\tilde x(t)+\D,\tilde u(t))-e(\tilde v(t))|^2_Q\nonumber\\
    -|e(\tilde x(t) + \Delta) - e(\tilde x(t))|^2_Q\nonumber\\
    +|g(\tilde x(t) +\Delta,\tilde u(t)) - \tilde y(t)|^2&\le\cl
    E_Q(\tilde z(t)) \label{mi48}
\end{flalign}
for all $\D$. Let $x(t)$ and $y(t)$ be 
defined by \rf{mi5},\rf{mi6} with
 $u(t)\equiv\tilde u(t)$ and
 $x(1)=\tilde x(1)$. 
 Substituting 
 $\D(t)= x(t)-\tilde x(t)$
 into \rf{mi48} yields
 \begin{align}
  &|e(x(t+1)) - e(\tilde x(t+1))|_Q^2 \nonumber\\
  & - |e(x(t)) - e(\tilde x(t))|_Q^2 + |y(t) - \tilde y(t)|^2
   \le {\cl E}_Q(\tilde z(t)).
\end{align}
Summing these inequalities over $t$  and noting:
$$|e(x(1))-e(\tilde x(1))|_Q^2 =0,\;\;|e(x(N+1))-e(\tilde
x(N+1))|_Q^2 \geq 0$$
yields \rf{mi49} .
\end{pf}

\subsection{Local RIE}
The local RIE for a DT model \rf{mi5},\rf{mi6} 
 is defined by:
\EQ{mi11}{
\cl E_Q^0(\tilde z)=
\sup_{\D}\left\{
|F\D+\e_x|^2_Q-|E\D|^2_Q+|G\D+\e_y|^2\right\},
}
and provides an upper bound for the linearized simulation error
$\bar{\cl E}^0$ according to the
following statement. 



\begin{thm}\label{thm:rieldt}
{\sl The inequality
\EQ{mi51}{  \bar{\cl E}^0\le\sum_{t=1}^N
\cl E_Q^0(\tilde z(t)),
}
holds
for every $Q=Q'>0$ and signal data \rf{mi44},\rf{mi45}.
}
\end{thm}
\begin{pf}
By the definition of ${\cl E}^0_Q$:
\EQ{eq:dtldiss}{|F\D+\e_x|^2_Q-|E\D|^2_Q+|G\D+\e_y|^2 \leq {\cl E}^0_Q(\tilde
z)}
holds for all $\Delta$.  Substituting $\Delta(t) = \tilde \Delta(t)$
defined by \eqref{eq:dtldiss}, with $\tilde \Delta(1) = 0$, we have:
$$|E(\tilde v(t))\tilde \D(t+1)|^2_Q-|E(\tilde
x(t))\D(t)|^2_Q+|G\tilde\D(t)+\e_y|^2 \leq {\cl E}^0_Q(\tilde
z).$$
Summing over $t$ yields \eqref{mi51}.
\end{pf}

Note that 
the supremum in \rf{mi11} is finite only when
the matrix
\EQ{mi68}{ R_{dt}=F'QF-E'QE+G'G}
is negative semidefinite. In applications, {\sl strict}
negative definiteness of the matrix \rf{mi68} 
is enforced, to be referred to as {\sl robustness} of the
corresponding supremum.

\subsection{RIE and Stability}
The following theorem shows that global
finiteness of the local RIE implies global stability
of the model \rf{mi5},\rf{mi6}.

\begin{thm}
{\sl 
Let continuously differentiable functions
$f,g,e$ and matrix $Q=Q'>0$ be 
such that $e$ has a smooth inverse $e^-$ (i.e. $e^-(e(x))=e(e^-(x))=x$
for all $x\in\RR^n$),
and $\cl E_Q^0(e^-(f(x,u)),x,u,g(x,u))$ is finite
for every $x\in\RR^n$, $u\in\RR^m$. 
Then system \rf{mi5},\rf{mi6} is
globally incrementally output $\ell_2$-stable. \label{thm:riestab}
}
\end{thm}
\begin{pf}
Let $(u,x,y)=(u_0,x_0,y_0)$ and 
$(u,x,y)=(u_1,x_1,u_1)$ be two solutions of
\rf{mi5},\rf{mi6} with $u_0=u_1=u$. For $\t\in[0,1]$
define $(x_*(\theta,t),y_*(\theta,t))$ 
as the solution
of \rf{mi5},\rf{mi6} with 
\[
x_*(\theta,1)=\t x_1(1)+(1-\t)x_0(1).\] 
Then $x_*(\t,t)$,
$y_*(\t,t)$  are continuously differentiable functions of $\t\in[0,1]$
for all integer $t\ge1$, and 
\[  y_*(0,t)=y_0(t),\ \  \ y_*(1,t)=y_1(t)\ \ \forall\ t\ge0.\]

Differentiating the identities
\begin{eqnarray*}  
e(x_*(\t,t+1))&=&f(x_*(\t,t),u(t)),\\
y_*(\t,t)&=&g(x_*(\t,t),u(t))
\end{eqnarray*}
with respect to $\t$ yields
\begin{eqnarray*}
E(x_*(\t,t+1))\frac{\partial x_*(\t,t)}{\partial\t}&=&
F(x_*(\t,t),u(t))\frac{\partial x_*(\t,t)}{\partial\t},\\
\frac{\partial y_*(\t,t)}{\partial\t}&=&
G(x_*(\t,t),u(t))\frac{\partial x_*(\t,t)}{\partial\t}.
\end{eqnarray*}
Since the finiteness of $\cl E_Q^0(e^-(f(x,u)),x,u,g(x,u))$
implies negative semidefiniteness of the quadratic
form
\[  \s(\D)=|F(x,u)\D|^2_Q-|E(x)\D|^2_Q+|G(x,u)\D|^2\]
for all $x,u$, we have
\begin{equation}
\label{eq:dissipation}
w(t)\le V(t)-V(t+1)
\end{equation}
for all $t\ge0$, where
\begin{eqnarray*}
w(t)&=&\int_0^1
\left|\frac{\partial y_*(\t,t)}{\partial\t}\right|^2d\t
\ge|y_*(0,t)-y_*(1,t)|^2,\\
 V(t)&=&\int_0^1
\left|\frac{\partial e(x_*(\t,t))}{\partial\t}\right|^2d\t.
\end{eqnarray*}
Summing \eqref{eq:dissipation}  over $t$ we find:
\[
\sum_{t=1}^N w(t)\le V(1) - V(N+1),\] and as $V(N+1) \geq 0$ the sum of
$w(t)$ is finite for all $N$. Since $w(t)\ge|y_0(t)-y_1(t)|^2$,
this proves incremental L2 output stability.
\end{pf}


\section{A Convex Upper Bound for Optimization}
The results of the previous section
suggest minimization (with respect to
$e,f,g,Q$) of the sum of RIE over the available
data points as an approach to system identification.
However, in general, the RIE functions are {\sl not}
convex with respect to $e$,
$f$, $g$ and
$Q$. In this section, we use the inequality
\EQ{mi61}{ -a'Qa\le \D'Q^{-1}\D-2\D'a,}
which, due to the identity
\[\D'Q^{-1}\D-2\D'a+a'Qa=
|a-Q^{-1}\D|^2_Q,\]
is valid for all $a,\D\in\RR^n$ and a real 
symmetric $n$-by-$n$ matrix
$Q$ such that $Q=Q'>0$, to derive a family of
upper bounds for the RIE functions. The upper bounds
will be jointly convex
with respect to $e$, $f$, $g$, and $P=Q^{-1}>0$.

 \subsection{Upper Bounds for Global RIE in Discrete Time}
Given a symmetric positive definite
$n$-by-$n$ matrix $Q$
and functions $e:\ \RR^n\mapsto\RR^n$, 
$f:\ \RR^n\times\RR^m\mapsto\RR^n$ let 
\begin{eqnarray*}
 \d_v   &=&f(\tilde x+\D,\tilde u)-e(\tilde v).
\end{eqnarray*}
Applying \rf{mi61} with $a=\d_e$, 
to the $-|\d_e|^2_Q$ term in the definition of $\cl E_Q(\tilde z)$
yields 
$\cl E_Q(\tilde z)\le\hat{\cl E}_Q(\tilde z)$ where
\EQ{mi63}{
\hat{\cl E}_Q(\tilde z)=\sup_{\D}\left\{|\d_v|_Q^2
+|\Delta|^2_P-2\Delta'\d_e
+|\d_y|^2\right\},
}
and $P=Q^{-1}$.
The function $\hat{\cl E}_Q(\tilde z)$ serves as an upper bound
for $\cl E_Q(\tilde z)$ that is jointly convex with respect to
$e$, $f$, $g$, and $P=Q^{-1}>0$.

\subsection{Upper Bounds for Local RIE in Discrete Time}
Given a symmetric positive definite
$n$-by-$n$ matrices $Q$
and functions $e:\ \RR^n\mapsto\RR^n$, 
$f:\ \RR^n\times\RR^m\mapsto\RR^n$ let
\begin{eqnarray*}  
\D_e&=&E(\tilde x)\D,\\
\D_v&=&F(\tilde x,\tilde u)\D+\e_x,\\
\D_y&=&G(\tilde x,\tilde u)\D+\e_y.
\end{eqnarray*}
Applying \rf{mi61} with $a=\D_e$, 
to the $-|\D_e|^2_Q$ term in the definition of $\cl E_Q(\tilde z)$
yields 
$\cl E_Q(\tilde z)\le\hat{\cl E}_Q^0(\tilde z)$ where
\EQ{mi23}{
\hat{\cl E}_Q^0(z)=\sup_{\D}\left\{|\D_v|_Q^2
+|\D|^2_P-2\D'{\D_e}
+|\D_y|^2\right\},
}
with $P=Q^{-1}$.
The function $\hat{\cl E}_Q^0(\tilde z)$ serves as an upper bound
for $\cl E_Q^0(\tilde z)$ that 
is jointly convex with respect to
$e$, $f$, $g$, and $P=Q^{-1}>0$.

\subsection{Well-Posedness of State Dynamics}
The well-posedness of state dynamics equation \rf{mi5}
is guaranteed when the function $e:\ \RR^n\mapsto\RR^n$
is a bijection. 
The well-posedness of \rf{mi5} is implied by robustness
of the supremum in the definition \rf{mi23} of the upper bound
$\hat{\cl E}_Q^0$ of the local RIE $\cl E_Q^0$, i.e. by
strict negative definiteness of the matrix:
\EQ{mi69}{\hat R_{dt}=F'QF+P-E'-E+G'G.}
Note that this is not guaranteed by the robustness of \eqref{mi11}.

\begin{thm}\label{thm:wellpos}
Let $e: \RR^n \mapsto \RR^n$ be a continuously differentiable function 
with a uniformly bounded Jacobian $E(x)$, satisfying:
\EQ{Ebnd}{E(x) + E(x)' \geq 2r_0I, \quad \forall x \in \RR^n}
for some fixed $r_0 > 0$.  Then $e$
is a bijection.  
\end{thm}
\begin{pf}

Consider the task of minimizing $|e(x)-z|^2$ with respect to
$x\in\RR^n$ for a given $z\in\RR^n$.
Since $E+E'\geq2r_0I$ implies 
\begin{eqnarray*}
 \frac{d}{d\theta}\D'[e(x+\D \theta)-e(x)]&=&\D'E(x+\D \theta)\D \\
&\ge & r_0|\D|^2,
\end{eqnarray*}
we have
\EQ{mi71}{ |e(x+\D)-e(x)|\ge r_0|\D|\ \ \forall\ x,\D,}
hence $|e(x)|\to\infty$ as $|x|\to\infty$, and the minimum
of $|e(x)-z|^2$ is achieved at some $x=x_0$. Then
the first order optimality condition
$(e(x_0)-z)'E(x_0)=0$ implies $e(x_0)=z$.
To show that the equation $e(x)=z$ has a unique
solution, use \rf{mi71}.
\end{pf}

When $e(x)$ is nonlinear one can solve for $\hat x$ such that $|\hat x
- x_0| < \epsilon$ (with $e(x_0) = z$) via the ellipsoid method, or
related techniques.  Given a guess $\hat x$, we know the true solution
lies in a sphere: $|e(\hat x) - z| \geq r_0|\hat x - x_0|$.
Further, we have a cutting plane oracle: $(\hat x - x_0) '(e(\hat x) -
z) \geq 0$.


\subsection{Coverage of Stable Linear Systems}
\label{sec:dtlinear}
Since we have produced an upper bound for 
the simulation error both through the 
introduction of $\cl E_Q(\tilde z)$ and 
$\hat{\cl E}_Q(\tilde z)$, it is desirable 
to check whether a basic class of systems will be recovered exactly.

Consider a linear system
\[
x(t+1) = Ax(t)+Bu(t),  \ y(t) = Cx(t)+Du(t)
\]
where $x\in \mathbb R^n$ and $u\in\mathbb R^m$. Define the ``data
matrices'' from an experiment of length $N$ to be $X:=[\tilde x(t_1), \ldots, \tilde x(t_N)]$, $U:=[\tilde u(t_1), \ldots, \tilde u(t_N)]$. Suppose we have fit a linear model
\[
Ex(t+1) = Fx(t)+Lu(t), \ y(t) = Gx(t)+Hu(t).
\]
We consider a linear system to have been recovered exactly by the model if $G=C, D=H, EB=L,$ and $EA=F$.

\

\begin{thm}\label{thm:linear}
For data generated from a stable DT linear system with zero noise, if the data matrix $[X', U']'$ is of rank at least $n+m$, then the linear system is recovered exactly and 
\begin{equation}\notag
\cl E_Q(\tilde z)= \hat{\cl E}_Q(\tilde z)=0.
\end{equation}
\end{thm}

\ 

Note that by construction for the case of a 
linear model $\cl E_Q = \cl E_Q^0$. In order to prove this theorem, we will use the following lemma:
\begin{lem}\label{lem:linear}
For any Schur matrix $A$ there exists $E, F$ and $Q>0$ such that $EA=F$ satisfying $M=M'<0$ where
\begin{equation}\label{eqn:parabola_EF}
M:=F'QF+Q^{-1}-E'-E+G'G.
\end{equation}
\end{lem}
\begin{pf}
Since $A$ is Schur, there exists a matrix $R>0$ 
such that $A'RA-R<-G'G$.  Let $E=R, F = RA, Q = R^{-1}$.
Substituting into (\ref{eqn:parabola_EF}) results in 
\begin{equation}
M=A'RA-R+G'G<0 \notag
\end{equation}
where the last inequality follows by construction of $R$.
\end{pf}

{\em Proof of Theorem \ref{thm:linear}.}
Using the choice of $E, F, Q$ in Lemma \ref{lem:linear}, since the data is noise free and $EA=F$ we have $\e_x=\e_y=0$. As a result, it follows from Lemma \ref{lem:linear} that $\hat{\cl E}_Q(\tilde z)$ is the supremum of a homogeneous negative-definite quadratic form in $\Delta$, hence has a value of zero. Similarly, with zero noise (\ref{mi11}) is the supremum of a homogeneous quadratic form in $\Delta$ and since $\cl E_Q(\tilde z)\le \hat {\cl E}_Q(\tilde z)=0$ hence $\cl E_Q(\tilde z)=0$. The rank condition on the data matrices ensures that if robust equation error is zero, then the true system is recovered.\hfill$\Box$\smallbreak




 \section{Continuous Time Results}
 For continuous time (CT) models, \rf{mi5}
is replaced by
\begin{eqnarray}
\frac{d}{dt} e(x(t))&=&f(x(t),u(t)),\label{mi7}
\end{eqnarray}
or, equivalently,
\begin{eqnarray*}
E(x(t))\dot x(t)&=&f(x(t),u(t)),
\end{eqnarray*}
where $E(x)$ is the Jacobian of $e(\cdot)$ at $x$.
Naturally, $E(x)$ is required to be non-singular for all $x$.
We consider the model 
\rf{mi6},\rf{mi7} is {\sl stable}
if the difference $y_1-y_0$ is
square integrable for every two solutions $(u,x,y)=(u_1,x_1,y_1)$
and $(u,x,y)=(u_0,x_0,y_0)$ of \rf{mi6},\rf{mi7} with
the same input $u_1=u_0=u$.

We expect to have input/state/output information in the form
$\tilde {\cl Z}=\{\tilde z(t_i)\}_{i=1}^N$, where
$\tilde z(t_i)=(\tilde v(t_i),\tilde x(t_i),\tilde u(t_i),\tilde
y(t_i))$. 
Here $\tilde x, \tilde u,\tilde y$ represent state, input and
output respectively, whereas $\tilde v \sim \dot {\tilde x}$.

For the purpose of theoretical analysis, we will assume that the
input/state/output information is available in the form of {\it
  signals}, that is functions: \EQ{eq:ctsig1}{ \tilde x, \tilde v:\
  {\cl T} \mapsto \R^n,\, \tilde u:\ {\cl T} \mapsto \R^m,\, \tilde
  y:\ {\cl T} \mapsto \R^k } such that \EQ{eq:ctsig2}{ \cl T = [0,T],
  \ \ \ \tilde v(t) = \frac{d}{dt} \tilde x(t), \forall\ t \in [0,T].
} In practice we have only sampled data, but for
theoretical convenience we assume $\tilde u(t)$ and $\tilde v(t)$
exist as suitably smooth functions (e.g. piecewise continuous)
interpolating the samples.

\subsection{Simulation Error}
 Given CT signal data
$\tilde {\cl Z}$,
and functions $e,f,g$, the {\sl simulation error}
associated with a model matching \rf{mi6},\rf{mi7} is defined
as 
\EQ{mi22}{ \bar{\cl E}=\int_0^T|\tilde y(t)- y(t)|^2dt,}
where $y$ is defined by \rf{mi6},\rf{mi7} with
$u(t)\equiv\tilde u(t)$ and
$x(0)=\tilde x(0)$.

\subsection{Linearized Simulation Error}
Similar to the DT case we examine a ``perturbed'' version of the
system equations:
\begin{eqnarray}
\frac{d}{dt} e(x_\t(t))&=&f(x_\t(t),u(t))-f_0(t),\label{eq:pertct1}\\
y_\t(t)&=&g(x_\t(t),u(t))-g_0(t).\label{eq:pertct2}
\end{eqnarray}
Here $\t \in [0,1]$, $f_0(t) = (1-\theta) \epsilon_x(\tilde z(t)),$ and $g_0(t) = (1-\theta)
\epsilon_y(\tilde z(t))$.  We examine the solution $(y_\t,x_\t)$ with $x_\t(0) = \tilde x(0)$ and $u(t) \equiv \tilde u(t)$.  For
the CT case, the equation error $\epsilon_x$ of \eqref{mi8} is replaced by:
\EQ{eq:eqerrct}{ \epsilon_x(\tilde z) = f(\tilde x,\tilde u)-E(\tilde x) \tilde v}
Note that for $\theta = 0$ we have $y_\theta = \tilde y$.

Via a linearized analysis similar to Section \ref{sec:linearizeddt} we have:
\EQ{mi38}{ \bar{\cl E}^0=
\int_0^T|G(\tilde x(t),\tilde u(t))\tilde\D(t)+
\e_y(\tilde z(t))|^2dt,}
where $\tilde\D(\cdot)$ is defined by
\EQ{mi39}{\frac{d}{dt}[ E(\tilde x(t))\tilde \D(t)]=
F(\tilde x(t),\tilde u(t))\tilde\D(t)+
\e_x(\tilde z(t)),}
with initial condition $\tilde\D(0)=0$.

\subsection{Global RIE in Continuous Time}
The global RIE error measure for a CT model \rf{mi6},\rf{mi7}
is similarly a function of $e,f,g$ and $Q = Q' > 0$, as well as a
single data-point $\tilde z$:
\EQ{mi14}{
\cl E_Q(\tilde z)=
\sup_{\D}\left\{2\d_e'Q[f(\tilde x+\D,\tilde u)-E(\tilde x)\tilde v]+
|\d_y|^2\right\}.
}

\begin{thm}\label{thm:riegct}
{\sl The inequality
\EQ{mi47}{  \bar{\cl E}\le\int_0^T
\cl E_Q(\tilde z(t))dt,
}
where $\tilde z(t)=(\tilde v(t),\tilde x(t),\tilde u(t),\tilde y(t))$,
holds
for every $Q=Q'>0$ and signal data \rf{eq:ctsig1},\rf{eq:ctsig2}.
}
\end{thm}
\begin{pf}
By the definition of $\cl E_Q(z)$ we have
\EQ{eq:ctgdiss}{ 2\d_e'Q[f(\tilde x+\D,\tilde u)-E(\tilde x)v]+|\D_y|^2\le\cl E_Q(\tilde z)}
for all $\D$. Let $(x,y)$ be
defined by \rf{mi6},\rf{mi7} with
$u(t)\equiv\tilde u(t)$ and
$x(0)=\tilde x(0)$. 
Substituting $\D=x(t)-\tilde x(t)$ into \rf{eq:ctgdiss} yields
\[ \frac{d|e(x(t))-e(\tilde x(t))|^2_Q}{dt}+
|\tilde y(t)- y(t)|^2\le \cl E_Q(\tilde z(t)).\]
Integrating this over the interval $t\in[0,T]$ yields \rf{mi47}.
\end{pf}

Theorem~\ref{thm:riegct} suggests minimization of the integral in
\rf{mi47} as an easier-to-handle alternative to minimization of
the simulation error.
In the case when system information comes in the {\sl sampled
data} format $\cl Z=\{z(t_i)\}_{i=1}^N$, the theorem
suggests minimization of the sum $\cl E_Q(z(t_i))$ with respect to
$Q=Q'>0$, $e$, $f$, $g$ as a system identification
algorithm.

\subsection{Local RIE in Continuous Time}
The local RIE error measure for a CT model \rf{mi6},\rf{mi7} 
is defined by
\EQ{mi12}{
\cl E_Q^0(z)=
\sup_{\D}\left\{
2(E\D)'Q(F\D+\e_x)+|G\D+\e_y|^2\right\}
}
and provides an upper bound for the linearized simulation error
$\bar{\cl E}^0$ according to the
following statement. 

\begin{thm}\label{thm:rielct}
{\sl The inequality
\EQ{mi50}{  \bar{\cl E}^0\le\int_0^T
\cl E_Q^0(\tilde z(t))dt,
}
holds
for every $Q=Q'>0$ and signal data \rf{eq:ctsig1},\rf{eq:ctsig2}.
}
\end{thm}

Note that  the supremum in \rf{mi12} is finite only when
the matrix \EQ{mi67}{ R_{ct}=E'QF+F'QE+G'G} is negative semidefinite. 

\subsection{RIE and Stability}
A similar statement to Theorem \ref{thm:riestab} is available in the CT case:
\begin{thm}
{\sl 
Let two times continuously differentiable functions
$e,f,g$ and matrix $Q=Q'>0$ be 
such that $E(x)$ is invertible
for all $x\in\RR^n$),
and $\cl E_Q^0(E(x)^{-1}f(x,u),x,u,g(x,u))$ is finite
for every $x\in\RR^n$, $u\in\RR^m$. 
Then system \rf{mi6},\rf{mi7} is
globally incrementally output ${\cl L}_2$-stable.
}
\end{thm}

\subsection{Upper Bounds for Continuous Time Global RIE}
Given a symmetric positive definite
$n$-by-$n$ matrix $Q$
and functions $e:\ \RR^n\mapsto\RR^n$, 
$f:\ \RR^n\times\RR^m\mapsto\RR^n$ let
\begin{eqnarray*}
\delta_e^+ &=& \d_e+f(\tilde x+\Delta,\tilde u) -E(\tilde x)\tilde v,\\
\delta_e^- &=& \d_e-f(\tilde x+\Delta,\tilde u) +E(\tilde x)\tilde v,
\end{eqnarray*}
where $E$ is the Jacobian of $e$.

Applying \rf{mi61} with $a=\d_e^-$, 
to the second term in the expression on the right side of
the identity
\[ 4\d_e'Q[f(\tilde x+\Delta,\tilde u) -E(\tilde x)\tilde v]=|\d_e^+|^2_Q-|\d_e^-|^2_Q\]
yields
$\cl E_Q(\tilde z)\le\hat{\cl E}_Q(\tilde z)$ where
\EQ{mi62}{
\hat{\cl E}_Q(\tilde z)=\sup_{\D}\left\{\frac{|\d_e^+|_Q^2
+|\Delta|^2_P}2-\Delta'\d_e^{-}+|\d_y|^2\right\},
}
that $P=Q^{-1}$.
The function $\hat{\cl E}_Q(\tilde z)$ serves as a CT upper bound
for $\cl E_Q(\tilde z)$ that is jointly convex with respect to
$e$, $f$, $g$, and $P=Q^{-1}>0$.

\subsection{Upper Bounds for Continuous Time Local RIE}
Given a symmetric positive definite
$n$-by-$n$ matrices $Q$
and functions $e:\ \RR^n\mapsto\RR^n$, 
$f:\ \RR^n\times\RR^m\mapsto\RR^n$ let
\begin{eqnarray*}
\D_e^+ &=& E(\tilde x)\D+F(\tilde x,\tilde u)\D+\e_x,\\
\D_e^- &=& E(\tilde x)\D-F(\tilde x,\tilde u)\D-\e_x,
\end{eqnarray*}
where $E,F,G$ are the Jacobians of $e,f,g$
with respect to $x$. 

Applying \rf{mi61} with $a=\D_e^-$, 
to the second term in the expression on the right side of
the identity
\[ 4(E\D)'Q[F\D+\e_x]=|\D_e^+|^2_Q-|\D_e^-|^2_Q\]
yields
$\cl E_Q(\tilde z)\le\hat{\cl E}_Q^0(\tilde z)$ where
\EQ{mi64}{
\hat{\cl E}_Q^0(z)=\sup_{\D}\left\{\frac{|\D_e^+|_Q^2
+|\D|^2_P}2-\D'\D_e^-+|\D_y|^2\right\},
}
with $P=Q^{-1}$. 
The function $\hat{\cl E}_Q^0(\tilde z)$ serves as a CT upper bound
for $\cl E_Q^0(\tilde z)$ that is 
jointly convex with respect to
$e$, $f$, $g$, and $P=Q^{-1}>0$.

\subsection{Well-Posedness of State Dynamics}
A CT model is well posed so long as $e$ from \rf{mi7} has a non-singular Jacobian $E=E(x)$
at every point $x\in\RR^n$.  Invertibility of the Jacobian at a given
point $x$ is guaranteed by {\sl robustness} of the supremum in the definition
\rf{mi12} of the local RIE $\cl E_Q^0$ (i.e. strict negative
definiteness of $R_{ct}$ in \eqref{mi67}).

\subsection{Recovery of Linear Systems}
A result similar to Theorem \ref{thm:linear} can also be shown in the
CT case based on the following lemma.
\begin{lem}\label{lem:ctlinear}
  For any Hurwitz matrix $A$ there exists $E,F$ and $Q = Q' > 0$ such
  that $F = EA$ and $M = M' < 0$ where:
\begin{flalign}\label{eqn:ctparabola_EF}
  M:=&(E+F)'Q(E+F) +Q^{-1}\notag\\
  &\quad -(E-F)'-(E-F)+2G'G. 
\end{flalign}
\end{lem}
\begin{pf}
  Since $A$ is Hurwitz, there exists an $R = R' > 0$ such that $A'R +  RA < -G'G$.  Take
  $E = (I-A)'R$,  $F = (I-A)'RA$, and  $Q = ((I-A)'R(I-A))^{-1}$
Note that as $A$ is Hurwitz, $I-A$ will be nonsingular.  Substituting
these choices into \eqref{eqn:ctparabola_EF} we have:
$$ M = 2A'R + 2RA  + 2G'G < 0$$
where the last inequality holds by the construction of $R$.
\end{pf}

We again consider ``data matrices'' $X:=[\tilde x(t_1), \ldots, \tilde x(t_N)]$, and $U:=[\tilde u(t_1), \ldots, \tilde u(t_N)]$.

\begin{thm}\label{thm:linear}
For data generated from a stable CT linear system with zero noise, if
the data matrix $[X', U']'$ is of rank at least $n+m$, then the linear system is recovered exactly and 
\begin{equation}\notag
\cl E_Q(\tilde z)= \hat{\cl E}_Q(\tilde z)=0.
\end{equation}
\end{thm}
\begin{pf}
  The proof is nearly identical to that of Theorem \eqref{thm:linear},
  using Lemma \eqref{lem:ctlinear} as necessary.
\end{pf}


\section{Implementation Details}
We now discuss practical considerations for data preparation and
minimization of the upper bounds using semidefinite programming.
\subsection{Approximating States}
\label{sec:state}
The RIE formulation assumes access to approximate state observations,
$\tilde x(t)$.  In most cases of interest, the full state of the
system is not directly measurable.  In practice, our solutions have
been motivated by the assumption that future output can be
approximated as a function of recent input-output history and future
input.  To summarize recent history, we have had success applying
linear filter banks, as is common in linear identification
(e.g. Laguerre filters \cite{Chou99}).

Even in fairly benign cases one expects the input-output histories to
live near a nonlinear submanifold of the space of possible histories.  As a
result, linear projection based methods may require excessive
dimensionality to approximate the state of the system.  Connections
between nonlinear dimensionality reduction and system identification
are being explored in the manifold learning community, such as
\cite{Rahimi06} and \cite{Ohlsson08}.

For CT identification estimating the rates of the
system, $v(t) =\frac{d}{dt} x(t)$, presents an additional
challenge. For true system outputs, this can be approached via
differentiation filters, or noncausal smoothing before numerical
differentiation.  Approximating additional states
through filter banks allows the rates of these variables to be
calculated analytically.



\subsection{Quality of Fit with Semidefinite Programs}
For any tuple of data, $\tilde z(t_i)$, the upper bound on the local RIE
is the supremum of a concave quadratic form in $\Delta$.  So long as
$e,f$ and $g$ are chosen to be linear in the decision variables, this upper
bound can be minimized by introducing an LMI for each data-point using
the Schur complement.  We introduce a slack variable $s_i$ for each
data-point: \EQ{eq:lcslack}{ s_i \geq \hat {\cal E}^{0}_Q(\tilde
z(t_i)),} which is a convex constraint and optimize for $\sum_i s_i
\rightarrow \min$.

Similarly, the upper bound on the global RIE is a function of $\Delta$
for fixed $\tilde z(t_i)$. If we take $e,f$ and $g$ to be
polynomials or rational functions with fixed denominators then the
upper bound will be a polynomial or rational function in $\Delta$.  As
a result, we can minimize this function by introducing a
sum-of-squares (SOS) constraint \cite{Parrilo00}. We again introduce a
slack variable $s_i$: \EQ{eq:glslack}{ s_i \geq \hat {\cal E}_Q
(\tilde z(t_i)),} and optimize for $\sum_i s_i \rightarrow \min$.  This
equation will be polynomial in $\Delta$ and quadratic in $n+1$ other
variables due to the Schur complement.  In most cases, replacing the
positivity constraint with a SOS constraint is another convex relaxation.

When fitting a linear (affine) model for \eqref{mi5},\eqref{mi6} or
\eqref{mi6},\eqref{mi7} it is interesting to note that $\hat {\cl E}
= \hat {\cl E}^0$ and further the SDP can be posed to grow
only with the dimension of the state, rather than the number of data
points.  For example, in the linear DT case one can compute the
supremum \eqref{mi23} (assuming it is finite) as a quadratic form in
the data:
\[
\hat {\cal E}^0(\tilde z) =
\left \|\begin{bmatrix} 0 \\ \epsilon_x \\ \epsilon_y \end{bmatrix}\right\|_{H^{-1}}^2
H = \begin{bmatrix} E + E - P  & F' & G' \\ F & P & 0 \\ G& 0& I\end{bmatrix}
\]
When minimizing the RIE over many data-points one can use the cyclic
property of trace to restate the problem in terms of the empirical
covariance matrix.  Using an eigenvalue decomposition of the
correlation matrix yields an
equivalent optimization problem with no more than $2n+m+k$ LMI constraints.

\subsection{Choice of Basis and Stability}
Global finiteness of the the upper bound $\hat {\cal E}^{0}_Q$
guarantees stability.  For a fixed $(x,u)$, boundedness can be verified
via an LMI.  Taking a polynomial or rational function basis for $e,f$ and
$g$, we can verify this LMI for all $(x,u)$ using a SOS constraint.
Global verification of the inequalities places some constraints on the
degrees of these polynomials.  For example, in DT the degree of $E(x)$ must be
able to be twice that of $F(x,u)$ for the inequality to hold globally.

In continuous time, we use the following parametrization to allow for
global stability verification:
\EQ{imp1}{e(x) = \frac{\bar e(x)}{q(x)}, \qquad
f(x) = \frac{\bar f(x,u)}{q(x)p(u)}. } Here $q(x): \RR^n \mapsto \RR$ is
a fixed polynomial of degree $2d_x$ in each $x_i$ such that $q(x) \geq
1$. Similarly $p(u): \RR^m \mapsto \RR$ is of degree $2d_u$ in each
$u_i$, and $p(u) \geq 1$. The numerators, $\bar f(x,u)$ and $\bar
e(x)$ are polynomials whose coefficients are decision variables.  Both
$\bar e(x)$ and $\bar f(x,u)$ are degree $2d_x+1$ in each $x_i$ and
$\bar f$ is of degree $2d_u$ in each $u_i$.

With these choices of degrees, it is possible for the convex
relaxation to be satisfied for all $(x,u)$.  The positivity of the
expression can be tested via a SOS decomposition.  In particular, we
choose $q(x)$ and $p(u)$ to be nearly constant over the range of the
observed data.  For example, we take:
\EQ{eq:pq}{q(x) = (1+\|x\|_2^2)^{d_x} \quad p(u) = (1+\|u\|_2^2)^{d_u}}
In general, centering and normalizing the data drastically improves
numerical properties of the method.  Here, rescaling the data such
that it lies in a unit ball around the origin makes this choice of $q$
and $p$ apply more generally.

When global stability is not required, care must be taken to ensure
that solutions to the implicit form equations still exist.  In
continuous time this is guaranteed if $E(x)$ is invertible for all
$x$, and similarly it is guaranteed if $e(x)$ is invertible in
discrete time.  Both of these constraints can be satisfied by
requiring $E(x) + E(x)' \geq 2r_0 I$, which can again be enforced
using a SOS constraint.


\section{Examples and Discussion}
\subsection{Stability and Noise}
When confronted with large measurement noise, we have observed that
RIE minimization produces models which are more stable (e.g. damped
for linear systems) than the system being fit.  This is most evident
in highly resonant, or nearly marginally stable systems.  In these
situations, we have had success minimizing equation error while simply
requiring the local RIE to be finite at the sample points. Mitigating
this effect is a focus of future work.
\subsection{Simulated DT Example}
We consider a second-order nonlinear simulated discrete time system:
\begin{flalign}
  \begin{bmatrix} 2v_1+ v_2^2 v_1+\frac{1}{3} v_1^5\\
    v_1+2v_2 + v_1^2 v_2+ \frac{1}{3} v_2^5\\
    \end{bmatrix}=&\: \begin{bmatrix}0.4 & -0.9 \\ 0.9 &
      0.4\end{bmatrix} x + \begin{bmatrix}u \\ 0\end{bmatrix},
\end{flalign}
where $x = [x_1(t)\: x_2(t)]'$ and $v_i(t) = x_i(t+1)$.  For training
we excite the system with a chirp: $\tilde u(t) =
4\sin(2\pi\frac{10}{500^2} t^2)$ for $t \in \{1,\ldots,500\}$.  We observe
$\tilde y(t) = \tilde x(t) = x(t) + w(t)$, where $w(t)$ is zero mean,
Gaussian i.i.d. measurement noise with covariance = $0.0025I$.

\begin{figure}
 \centering
 \includegraphics[width=\columnwidth]{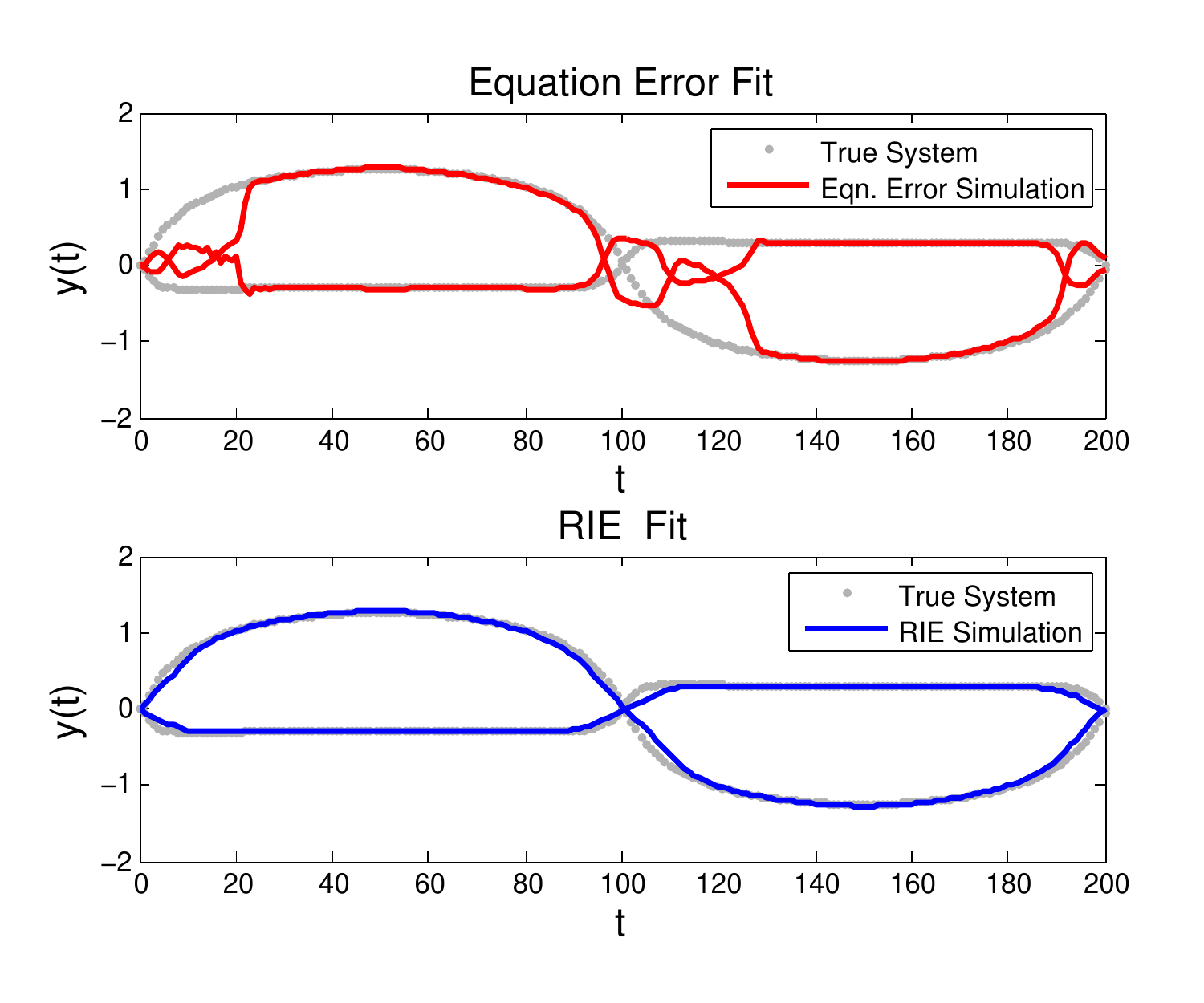}
 \caption{A comparison of equation error minimization and local RIE
   minimization on a simulated, second-order, nonlinear discrete time
   system. The true system response to validation input is compared to
 an equation error fit (top) and local RIE fit (bottom).  The system
 is not within the model class being searched over.}
 \label{fig:dt2_sim_fit}
 \end{figure}

We fit a model \eqref{mi5},\eqref{mi6} with
$g(\cdot,\cdot)$ fixed {\it a priori} to be $g(x,u) = x$.  We choose
$e(\cdot)$ to be cubic, and $f(\cdot,\cdot)$ to be a linear
combination of $u$ and the monomials up to total degree 7 in $x_i$.
With these choices the true system is outside the model class.
We compare minimizing the local RIE and minimizing equation
error.  In both cases, we restrict $E + E' > 2I$ to remove the scale
invariance of the problem.
Figure \eqref{fig:dt2_sim_fit} presents the response of the true system
and models for the input $u_{\textrm{test}}(t) =
4\sin(2\pi\frac{1}{200}t)$ over $t \in \{1,\ldots,200\}$.

\subsection{Modeling of Post-spike Dynamics in Live Neurons}
Our second example is drawn from the task of identifying the response
of the membrane potential of a live neuron. 
Details of the experimental procedure are given in the appendix. In particular, we are interested in
identifying the dynamics of the neuron immediately following an action
potential.  

We excite the neuron with $27$ separate multisine input currents.
The excitation is applied via a zero-order hold.  The response is the
sampled membrane potential of the neuron, $\tilde y(t)$. Both
measurement and control have a sampling rate of 10 kHz.  This data set
consists of 22 spikes which were separated into equal size training
and testing sets.

To achieve a 3rd order CT fit of the system, we pass the observed
output voltage, $\tilde y(t)$, through a filter bank determined by the first two
Laguerre functions with a pole at 300 radians per second \cite{Chou99}.  The
original voltage and the output of this filter bank give us $\tilde x(t)
\in \RR^3$. To compute $\tilde v(t)$ we apply a noncasual regularized
smoothing to the observed output and differentiate numerically.
For our model structure we choose $e,f$ polynomial in each $x_i$ (degree 4) and $f$ affine in $u$.
As our observation is a state, we fix our model's $g(x,u)$ to be the membrane potential.

As the response is nearly periodic, we avoid repetitive data by picking approximately 500 data points uniformly spread throughout the $(\tilde x, \tilde v, \tilde u)$ space.  We minimize $\hat {\cal E}^{0}_Q$.  
For comparison, we also fit a model of the same structure minimizing the equation error, $\sum_i |\epsilon_x(\tilde z(t_i))|^2 \rightarrow \min$.  In both cases, we insist on an invertible Jacobian $E(x)$ by requiring $E(x) + E(x)' \geq 10^{-3} I$ with $\delta = 1e-3$. 

\begin{figure} \centering \includegraphics[width=\columnwidth]{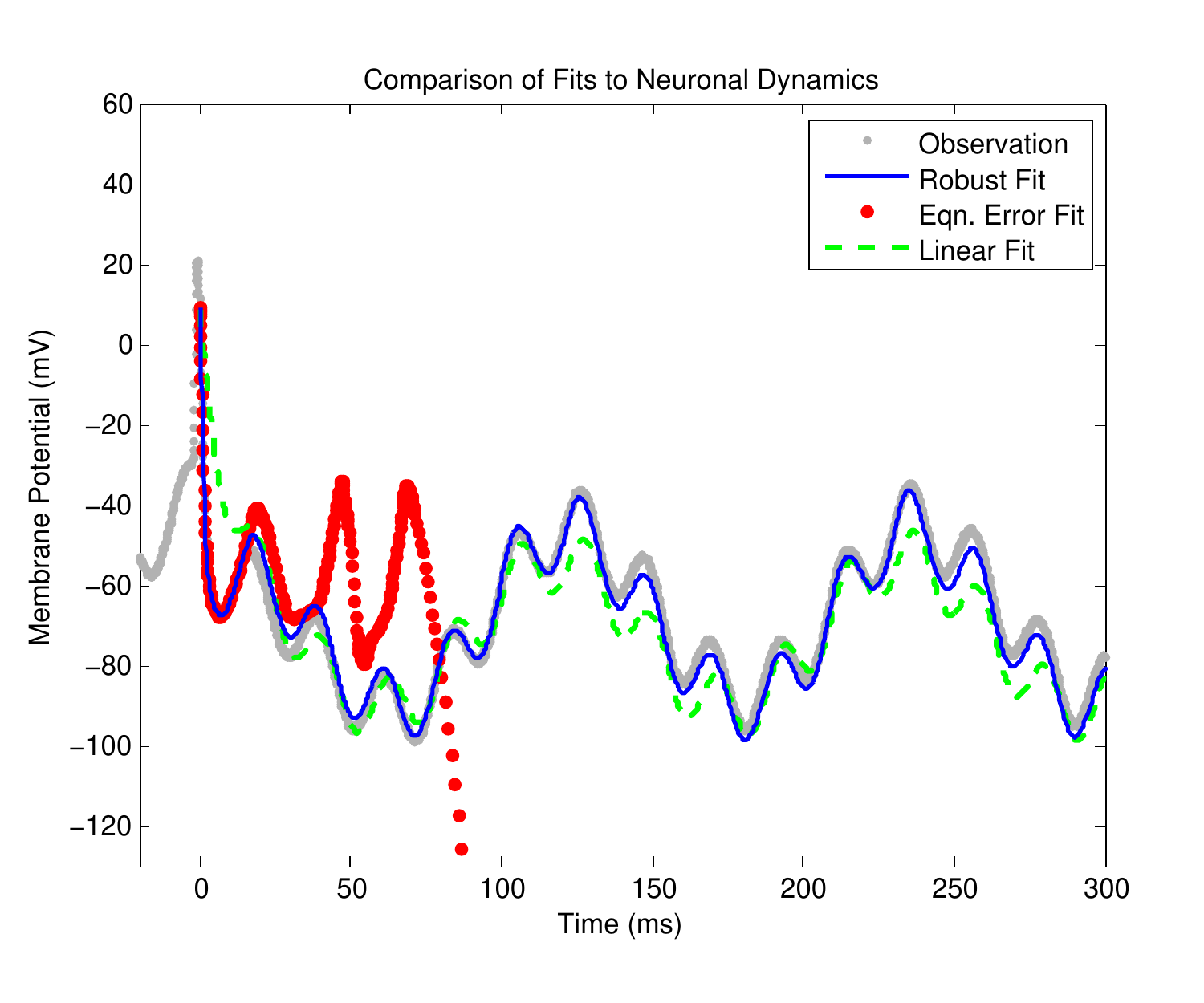} \caption{
We compare several fits of the post-spike dynamics of a live neuronal
cell on validation data.  The ``Robust Fit'' corresponds to minimizing
the local RIE, and is compared to both linear and nonlinear fits
minimizing equation error.  By $t= 100$, the nonlinear equation error
fit has diverged.   The linear fit
does not capture the steep descent at $t=0$, nor does it replicate the
long term behavior.}  \label{fig:neuron-fit} \end{figure}

Figure \ref{fig:neuron-fit} plots a neuronal response from the test
set and the result of simulating the models from the same initial 
conditions. Also included is a first order DT model fit using equation
error (CT and higher order linear equation error fits led to unstable models).

\appendix
\subsection{Live Neuron Experimental Procedure}

Primary rat hippocampal cultures were prepared from P1 rat pups, in
accordance with the MIT Committee on Animal Care policies for the
humane treatment of animals. Dissection and dissociation of rat hippocampi
were performed in a similar fashion to \cite{Hagler:2001p112}. Dissociated
neurons were plated at a density of 200K cells/mL on 12 mm round glass
coverslips coated with 0.5 mg/mL rat tail collagen I (BD Biosciences)
and 4 $\mu$g/mL poly-D-lysine (Sigma) in 24-well plates. After 2
days, 20 $\mu$M Ara-C (Sigma) was added to prevent further growth
of glia. 

Cultures were used for patch clamp recording after 14 days in vitro.
Patch recording solutions were previously described in \cite{Bi:1998p116}.
Glass pipette electrode resistance ranged from 2-4 M\textgreek{W}.
Recordings were established by forming a G\textgreek{W} seal between
the tip of the pipette and the neuron membrane. Perforation of the
neuron membrane by amphotericin-B (300 $\mu$g/mL) typically occurred
within 5 minutes, with resulting access resistance in the range of
10-20 M\textgreek{W}. Recordings with leak currents smaller than -100
pA were selected for analysis. Leak current was measured as the current
required to voltage clamp the neuron at -70 mV. Synaptic activity
was blocked with the addition of 10 $\mu$M CNQX, 100 $\mu$M APV,
and 10 $\mu$M bicuculline to the bath saline. 
%
Holding current was
applied as necessary to compensate for leak current.

\bibliographystyle{IEEEtran}
\bibliography{sysid}

\end{document}